\documentclass[a4paper,notitlepage,11pt]{article}%
\usepackage{amssymb}
\usepackage{graphicx}
\usepackage{amsmath}
\usepackage{amsthm}
\usepackage{amsfonts}
\usepackage{amssymb}%
\providecommand{\U}[1]{\protect\rule{.1in}{.1in}}
\theoremstyle{plain}
\newtheorem{theorem}{Theorem}[section]

\newtheorem{corollary}[theorem]{Corollary}

\newtheorem{lemma}[theorem]{Lemma}

\theoremstyle{definition}

\newtheorem{remark}[theorem]{Remark}
\numberwithin{equation}{section}
\numberwithin{theorem}{section}
\ifx\pdfoutput\relax\let\pdfoutput=\undefined\fi
\newcount\msipdfoutput
\ifx\pdfoutput\undefined\else
\ifcase\pdfoutput\else
\ifx\paperwidth\undefined\else
\ifdim\paperheight=0pt\relax\else\pdfpageheight\paperheight\fi
\ifdim\paperwidth=0pt\relax\else\pdfpagewidth\paperwidth\fi
\fi\fi\fi
\begin{document}

\title{On Dirichlet problems with singular nonlinearity of indefinite sign
\thanks{2000 \textit{Mathematics Subject Clasification}. 35J75, 35J25, 35B09,
34B16.}  \thanks{Partially supported by
Secyt-UNC. } }
\author{T. Godoy, U. Kaufmann, \thanks{\textit{E-mail addresses. }godoy@mate.uncor.edu
(T. Godoy), kaufmann@mate.uncor.edu (U. Kaufmann, Corresponding Author).}
\and \noindent\\{\small FaMAF, Universidad Nacional de C\'{o}rdoba, (5000) C\'{o}rdoba,
Argentina}}
\maketitle

\begin{abstract}
Let $\Omega$ be a smooth bounded domain in $\mathbb{R}^{N}$, $N\geq1$, let
$K$, $M$ be two nonnegative functions and let $\alpha,\gamma>0$. We study
existence and nonexistence of positive solutions for singular problems of the
form $-\Delta u=K\left(  x\right)  u^{-\alpha}-\lambda M\left(  x\right)
u^{-\gamma}$ in $\Omega$, $u=0$ on $\partial\Omega$, where $\lambda>0$ is a
real parameter. We mention that as a particular case our results apply to
problems of the form $-\Delta u=m\left(  x\right)  u^{-\gamma}$ in $\Omega$,
$u=0$ on $\partial\Omega$, where $m$ is allowed to change sign in $\Omega$.

\end{abstract}

\section{Introduction}

Let $\Omega$ be a smooth bounded domain in $\mathbb{R}^{N}$, $N\geq1$, and let
$0\leq K,M\in L^{p}\left(  \Omega\right)  $ for some $p\geq2$ if $N=1$ and
$p>N$ when $N\geq2$. Our aim in this article is to consider existence and
nonexistence of solutions for singular problems of the form%
\begin{equation}
\left\{
\begin{array}
[c]{ll}%
-\Delta u=K\left(  x\right)  u^{-\alpha}-\lambda M\left(  x\right)
u^{-\gamma} & \text{in }\Omega\\
u>0 & \text{in }\Omega\\
u=0 & \text{on }\partial\Omega,
\end{array}
\right.  \label{prob}%
\end{equation}
where $\alpha,\gamma>0$ and $\lambda>0$ is a real parameter.

Concerning the results in this paper, in Section 3 we shall study (\ref{prob})
when $N=1$. More precisely, in Theorem \ref{nuevo} we shall give sufficient
conditions for the existence of solutions in the case $\alpha=\gamma$, and as
a consequence we shall derive existence results when $\alpha>\gamma$ in
Corollary \ref{coro}. A further result without any relation between $\alpha$
and $\gamma$ is presented in Theorem \ref{algo}, while in Theorem
\ref{necesaria} we prove necessary conditions on $\lambda$, $K$ and $M$ in
order to have solutions for (\ref{prob}) (see also Remark \ref{rem}, where we
also characterize the set of $\lambda^{\prime}$s such that (\ref{prob}) admits
a solution). We remark that as a particular case (taking $K:=m^{+}$,
$M:=m^{-}$, $\alpha=\gamma$ and $\lambda=1$) we are able to deal with problems
of the form%
\begin{equation}
\left\{
\begin{array}
[c]{ll}%
-\Delta u=m\left(  x\right)  u^{-\gamma} & \text{in }\Omega\\
u>0 & \text{in }\Omega\\
u=0 & \text{on }\partial\Omega,
\end{array}
\right.  \label{sipi}%
\end{equation}
where $m$ may change in $\Omega$ (see Corollary \ref{golazo}).

In Section 4 we shall consider the case $N\geq2$. We point out that although
all the results obtained for the one-dimensional problem have their
counterpart in $N$ dimensions, the conditions derived in Section 3 are sharper
or provide more information than the ones available when $N>1$. Sufficient
conditions are stated in Theorems \ref{bola} and \ref{algo n} and Corollary
\ref{coro n}, and in Theorem \ref{necesaria n} some necessary conditions are
presented. The particular case (\ref{sipi}) is covered in Corollary
\ref{golazo n}. Let us mention that these theorems shall be obtained employing
some estimates on linear problems combined with the sub and supersolution
method (see Theorem \ref{subsup} below) together with Schauder's fixed point theorem.

In order to relate our results to others already existing, let us first
mention that problems like (\ref{sipi}) have been extensively studied both
when $m>0$ (see e.g. \cite{tartar}, \cite{lazer}, \cite{royal} and its
references) and also if $m$ is nonnegative (e.g. \cite{gomes}, \cite{pino},
\cite{lair}), but to our knowledge no theorems are known when $m$ is allowed
to change sign in $\Omega$. In fact, even when (\ref{sipi}) is sublinear (that
is, $\gamma\in\left(  -1,0\right)  $) and one-dimensional, these kind of
problems are quite intriguing and involved and, as far as we know, only
recently existence of (strictly) positive solutions has started being studied
in detail when $m$ changes sign in $\Omega$ (see \cite{nodea}, \cite{ejde},
\cite{junows} and \cite{jesusultimo}; and \cite{plap} for the $p$-laplacian).
Let us also mention that nonnegative solutions of these semilinear problems
have been studied carefully in \cite{bandle}.

On the other side, problem (\ref{prob}) with $\alpha=0$ and $M\equiv1$ was
treated in for instance \cite{morel}, \cite{choi}, while in \cite{royal} it
was considered assuming that $K>0$ (under stronger regularity assumptions),
but we could not find any results in the literature in the case where $K$
vanishes in a subset of $\Omega$. We refer the reader to the nice review
papers \cite{handd} and \cite{hand} for further references, applications and
historical remarks concerning these types of singular elliptic problems, and
for specifically one-dimensional singular problems we refer to the book
\cite{ravi} and the references therein.

\section{Preliminaries and auxiliary results}

Let $\Omega\subset\mathbb{R}^{N}$ be a smooth bounded domain. For $h\in
L^{q}(\Omega)$, $q>1$, let $\phi\in W^{2,q}(\Omega)\cap W_{0}^{1,q}(\Omega)$
be the unique solution of
\begin{equation}
\left\{
\begin{array}
[c]{ll}%
-\Delta\phi=h\left(  x\right)  & \text{in }\Omega\\
\phi=0 & \text{on }\partial\Omega,
\end{array}
\right.  \label{lineal}%
\end{equation}
and let us write $\mathcal{S}:L^{q}\Omega)\rightarrow W^{2,q}(\Omega)$ for the
solution operator of (\ref{lineal}).

The two following lemmas provide some useful lower and upper bounds for
$\mathcal{S}\left(  h\right)  $ when $h$ is nonnegative. We set
\[
\delta_{\Omega}\left(  x\right)  :=dist\left(  x,\partial\Omega\right)
\text{.}%
\]

\begin{lemma}
\label{tomas}Suppose $\Omega:=\left(  a,b\right)  $, and let $0\leq h\in
L^{q}(\Omega)$ with $q>1$. Then in $\overline{\Omega}$ it holds that%
\begin{gather}
\mathcal{S}\left(  h\right)  \geq\frac{1}{b-a}\left(  \int_{a}^{b}h\left(
t\right)  \delta_{\Omega}\left(  t\right)  dt\right)  \delta_{\Omega
}:=\underline{c}\delta_{\Omega}\qquad\text{and}\label{l1}\\
\mathcal{S}\left(  h\right)  \leq\frac{1}{b-a}\max\left(  \int_{a}^{b}\left(
t-a\right)  h\left(  t\right)  dt,\int_{a}^{b}\left(  b-t\right)  h\left(
t\right)  dt\right)  \delta_{\Omega}:=\overline{c}\delta_{\Omega}. \label{l11}%
\end{gather}
Moreover, the inequalities (\ref{l1}) and (\ref{l11}) are sharp in the sense
that if $c_{1}\delta_{\Omega}\leq\mathcal{S}\left(  h\right)  \leq c_{2}%
\delta_{\Omega}$ in $\Omega$ for some $c_{1},c_{2}>0$, then $c_{1}%
\leq\underline{c}$ and $\overline{c}\leq c_{2}$.
\end{lemma}

\textit{Proof}. Let $\phi:=\mathcal{S}\left(  h\right)  $. It is easy to
verify that (even if $h$ changes sign in $\Omega$)
\begin{gather}
\phi\left(  x\right)  =\frac{x-a}{b-a}\int_{a}^{b}\int_{a}^{y}h\left(
t\right)  dtdy-\int_{a}^{x}\int_{a}^{y}h\left(  t\right)  dtdy=\nonumber\\
\frac{x-a}{b-a}\int_{a}^{b}\left(  b-t\right)  h\left(  t\right)  dt-\int
_{a}^{x}\left(  x-t\right)  h\left(  t\right)  dt. \label{fifi}%
\end{gather}
Also, if $t_{1},t_{2}\in\Omega$ with $t_{1}<t_{2}$ we may integrate over
$\left(  t_{1},t_{2}\right)  $ (see e.g. \cite{brezislibro}, Theorem 8.2) to
obtain that%
\[
\phi^{\prime}\left(  t_{1}\right)  -\phi^{\prime}\left(  t_{2}\right)
=-\int_{t_{1}}^{t_{2}}\phi^{\prime\prime}\left(  t\right)  dt=\int_{t_{1}%
}^{t_{2}}h\left(  t\right)  dt\geq0
\]
and therefore we find that $\phi$ is concave in $\Omega$.

Since $\delta_{\Omega}\left(  x\right)  =\min\left(  x-a,b-x\right)  $, using
(\ref{fifi}) we now observe that
\begin{gather*}
\phi\left(  \frac{a+b}{2}\right)  =\frac{1}{2}\int_{a}^{b}\left(  b-t\right)
h\left(  t\right)  dt-\int_{a}^{\frac{a+b}{2}}\left(  \frac{a+b}{2}-t\right)
h\left(  t\right)  dt=\\
=\frac{1}{2}\left(  \int_{a}^{\frac{a+b}{2}}\left(  t-a\right)  h\left(
t\right)  dt+\int_{\frac{a+b}{2}}^{b}\left(  b-t\right)  h\left(  t\right)
dt\right)  =\underline{c}\delta_{\Omega}\left(  \frac{a+b}{2}\right)
\end{gather*}
and hence the concavity of $\phi$ yields that $\phi\geq\underline{c}%
\delta_{\Omega}$ in $\overline{\Omega}$.

On the other hand, from (\ref{fifi}) we also get that
\[
\phi^{\prime}\left(  a\right)  =\int_{a}^{b}\frac{b-t}{b-a}h\left(  t\right)
dt,\qquad\phi^{\prime}\left(  b\right)  =\int_{a}^{b}\left(  \frac{b-t}%
{b-a}-1\right)  h\left(  t\right)  dt.
\]
Furthermore, if $h\not \equiv 0$, since $h\geq0$ it holds that $\phi^{\prime
}\left(  b\right)  <0<\phi^{\prime}\left(  a\right)  $ and thus recalling the
definition of $\overline{c}$ again by the concavity of $\phi$ we deduce that
$\phi\leq\overline{c}\delta_{\Omega}$ in $\overline{\Omega}$.

Since the final assertion of the lemma is clearly true this concludes the
proof. $\blacksquare$

\begin{lemma}
\label{morel}Let $0\leq h\in L^{q}(\Omega)$ with $q>N$. Then in $\overline
{\Omega}$ it holds that%
\begin{equation}
c_{\Omega}\left(  \int_{\Omega}h\left(  x\right)  \delta_{\Omega}\left(
x\right)  dx\right)  \delta_{\Omega}\leq\mathcal{S}\left(  h\right)  \leq
C_{\Omega}\left\Vert h\right\Vert _{L^{q}\left(  \Omega\right)  }%
\delta_{\Omega} \label{m1}%
\end{equation}
for some $c_{\Omega},C_{\Omega}>0$ depending only on $\Omega$.
\end{lemma}

\textit{Proof}. The first inequality in (\ref{m1}) appeared first in an
unpublished work by Morel and Oswald, and there is a nice proof in the paper
of Brezis and Cabr\'{e}, \cite{cabre}, Lemma 3.2. On the other side, set
$\phi:=\mathcal{S}\left(  h\right)  $, let $x\in\Omega$ and pick $y\in
\partial\Omega$ such that $\left\vert y-x\right\vert =\delta_{\Omega}\left(
x\right)  $. Since by the Sobolev imbedding theorems $\mathcal{S}:L^{q}%
\Omega)\hookrightarrow C^{1}(\overline{\Omega})$ and hence $\left\Vert
\left\vert \nabla\phi\right\vert \right\Vert _{L^{\infty}\left(
\Omega\right)  }\leq C_{\Omega}\left\Vert h\right\Vert _{L^{q}\left(
\Omega\right)  }$ for some $C_{\Omega}>0$ depending only on $\Omega$, from the
mean value theorem we find that%
\[
\left\vert \phi\left(  x\right)  \right\vert =\left\vert \phi\left(  x\right)
-\phi\left(  y\right)  \right\vert \leq\left\Vert \left\vert \nabla
\phi\right\vert \right\Vert _{L^{\infty}\left(  \Omega\right)  }\delta
_{\Omega}\left(  x\right)  \leq C_{\Omega}\left\Vert h\right\Vert
_{L^{q}\left(  \Omega\right)  }\delta_{\Omega}\left(  x\right)
\]
which in turn gives the second inequality in (\ref{m1}). $\blacksquare$

\qquad

Let $f:\Omega\times\left(  0,\infty\right)  \rightarrow\mathbb{R}$ be a
Carath\'{e}odory function (that is, $f\left(  \cdot,\xi\right)  $ is
measurable for all $\xi\in\left(  0,\infty\right)  $ and $f\left(
x,\cdot\right)  $ is continuous for $a.e.$ $x\in\Omega$). We consider now
singular problems of the form
\begin{equation}
\left\{
\begin{array}
[c]{ll}%
-\Delta u=f\left(  x,u\right)  & \text{in }\Omega\\
u=0 & \text{on }\partial\Omega
\end{array}
\right.  \label{sing}%
\end{equation}
in a suitable sense. We say that $v\in W_{loc}^{1,2}\left(  \Omega\right)
\cap C\left(  \overline{\Omega}\right)  $ is a \textit{subsolution} (in the
sense of distributions) of (\ref{sing}) if $v>0$ in $\Omega$, $v=0$ on
$\partial\Omega$, and
\[
\int_{\Omega}\left\langle \nabla v,\nabla\varphi\right\rangle \leq\int
_{\Omega}f\left(  x,v\right)  \varphi\qquad\text{for all }0\leq\varphi\in
C_{c}^{\infty}\left(  \Omega\right)  .
\]
Similarly, $w\in W_{loc}^{1,2}\left(  \Omega\right)  \cap C\left(
\overline{\Omega}\right)  $ is a \textit{supersolution} of (\ref{sing}) if
$w>0$ in $\Omega$, $w=0$ on $\partial\Omega$, and
\[
\int_{\Omega}\left\langle \nabla w,\nabla\varphi\right\rangle \geq\int
_{\Omega}f\left(  x,w\right)  \varphi\qquad\text{for all }0\leq\varphi\in
C_{c}^{\infty}\left(  \Omega\right)  .
\]
For the sake of completeness we state the following existence result (for a
proof, see \cite{loc}, Theorem 4.1).

\begin{theorem}
\label{subsup}Suppose there exist $v,w\in C^{1}\left(  \Omega\right)  $ sub
and supersolutions respectively of (\ref{sing}), satisfying that $v\leq w$ in
$\Omega$. Suppose also that there exists $g\in L_{loc}^{2}\left(
\Omega\right)  $ such that $\left\vert f\left(  x,\xi\right)  \right\vert \leq
g\left(  x\right)  $ for $a.e.$ $x\in\Omega$ and all $\xi\in\left[  v\left(
x\right)  ,w\left(  x\right)  \right]  $. Then there exists $u\in C^{1}\left(
\Omega\right)  \cap C\left(  \overline{\Omega}\right)  $ solution (in the
sense of distributions) of (\ref{sing}) with $v\leq u\leq w$, that is,
\[
\int_{\Omega}\left\langle \nabla u,\nabla\varphi\right\rangle =\int_{\Omega
}f\left(  x,u\right)  \varphi\qquad\text{for all }\varphi\in C_{c}^{\infty
}\left(  \Omega\right)  .
\]

\end{theorem}

\begin{remark}
\label{supersol}When $0\leq K\in L^{p}\left(  \Omega\right)  $ with
$K\not \equiv 0$ and $p>N$ one can readily check that (\ref{prob}) admits
arbitrarily large supersolutions. Indeed, let $\psi:=\mathcal{S}\left(
K\right)  >0$ and define also $\sigma\in\left(  0,1\right)  $ by
$\sigma:=1/\left(  1+\alpha\right)  $. We have that $\psi^{\sigma}\in
W_{loc}^{2,p}\left(  \Omega\right)  \cap C\left(  \overline{\Omega}\right)  $,
$\psi^{\sigma}=0$ on $\partial\Omega$, and a simple computation shows that for
all $\lambda>0$ and every $c\geq\sigma^{-1/\left(  1+\alpha\right)  }$ it
holds that
\begin{gather*}
-\Delta\left(  c\psi^{\sigma}\right)  =-c\sigma\psi^{\sigma-1}\Delta
\psi-c\sigma\left(  \sigma-1\right)  \psi^{\sigma-2}\left\vert \nabla
\psi\right\vert ^{2}\geq\\
-c\sigma\psi^{\sigma-1}\Delta\psi=c\sigma\psi^{\sigma-1}K\left(  x\right)
\geq K\left(  x\right)  \left(  c\psi^{\sigma}\right)  ^{-\alpha}\geq\\
K\left(  x\right)  \left(  c\psi^{\sigma}\right)  ^{-\alpha}-\lambda M\left(
x\right)  \left(  c\psi^{\sigma}\right)  ^{-\gamma}\qquad\text{in }%
\Omega^{\prime}%
\end{gather*}
for every $\Omega^{\prime}\subset\subset\Omega$, and therefore $c\psi^{\sigma
}$ is a supersolution of (\ref{prob}). $\blacksquare$
\end{remark}

\section{The one-dimensional case}

In this section we shall assume that
\[
\Omega:=\left(  a,b\right)  \qquad\text{for some }-\infty<a<b<\infty.
\]
For $1<p\leq\infty$, we define as usual $p^{\prime}$ by $1/p+1/p^{\prime}=1$
if $p$ is finite and $p^{\prime}=1$ if $p=\infty$. For $0\leq h\in
L^{p}\left(  \Omega\right)  $ we also introduce the following notation:%
\begin{equation}
h_{a}\left(  x\right)  :=\left(  x-a\right)  h\left(  x\right)  ,\qquad
h_{b}\left(  x\right)  :=\left(  b-x\right)  h\left(  x\right)  . \label{om}%
\end{equation}
Let us note that $\left\Vert h_{a}\right\Vert _{L^{p}\left(  \Omega\right)
}\leq\left(  b-a\right)  \left\Vert h\right\Vert _{L^{p}\left(  \Omega\right)
}$ if $p$ is finite and that $\left\Vert h_{a}\right\Vert _{L^{\infty}\left(
\Omega\right)  }=\left(  b-a\right)  \left\Vert h\right\Vert _{L^{\infty
}\left(  \Omega\right)  }$, and analogously for $h_{b}$. We also set%
\[
P^{\circ}:=\text{interior of the positive cone of }C^{1}\left(  \overline
{\Omega}\right)  .
\]

\begin{theorem}
\label{nuevo}Let $0\leq K,M\in L^{2}\left(  \Omega\right)  $ with
$K\not \equiv 0$. Assume $\alpha=\gamma$ and let $M_{a}$ and $M_{b}$ be given
by (\ref{om}).\newline(i) If $M\in L^{p}\left(  \Omega\right)  $ with $p\geq
2$, $\gamma\in\left(  0,\left(  p-1\right)  /p\right)  $ and%
\begin{gather}
\max\left(  \left\Vert M_{a}\right\Vert _{L^{p}\left(  \Omega\right)
},\left\Vert M_{b}\right\Vert _{L^{p}\left(  \Omega\right)  }\right)  \leq
c_{\gamma,p,a,b}\frac{\left(  \int_{a}^{b}K(t)\delta_{\Omega}(t)dt\right)
^{1+\gamma}}{\left(  \int_{a}^{b}K(t)dt\right)  ^{\gamma}},\text{\qquad
where}\label{m2}\\
c_{\gamma,p,a,b}:=\frac{\gamma^{\gamma}}{\left(  \gamma+1\right)  ^{\gamma+1}%
}\frac{\left(  1-\gamma p^{\prime}\right)  ^{1/p^{\prime}}}{\left(
b-a\right)  ^{\gamma+1/p^{\prime}}},\nonumber
\end{gather}
then (\ref{prob}) has a solution $u\in C^{1}\left(  \Omega\right)  \cap
C\left(  \overline{\Omega}\right)  $ for all $\lambda\leq1$; and $u\in
W^{2,q}\left(  \Omega\right)  \cap P^{\circ}$, $q>1$, whenever $K\delta
_{\Omega}^{-\gamma}\in L^{r}\left(  \Omega\right)  $ with $r>1$.\newline(ii)
If%
\begin{equation}
\max\left(  \int_{a}^{b}M_{a}(t)dt,\int_{a}^{b}M_{b}(t)dt\right)  <\int
_{a}^{b}K\left(  t\right)  \delta_{\Omega}\left(  t\right)  dt, \label{hip}%
\end{equation}
then there exists $\gamma_{0}>0$ such that the problem (\ref{prob}) has a
solution $u\in W^{2,q}\left(  \Omega\right)  \cap P^{\circ}$, $q>1$, for all
$\gamma\in\left(  0,\gamma_{0}\right]  $ and all $\lambda\leq1$.
\end{theorem}

\textit{Proof}. Taking into account Theorem \ref{subsup} and Remark
\ref{supersol} we note that it is enough to build a subsolution for
(\ref{prob}). Furthermore, clearly any solution of (\ref{prob}) is a
subsolution of (\ref{prob}) with $\underline{\lambda}$ in place of $\lambda$
whenever $\underline{\lambda}\leq\lambda$ and so in order to prove the theorem
we may assume that $\lambda=1$.

Let us prove (i). Since (\ref{prob}) with $\alpha=\gamma$ is homogeneous (i.e.
it has a solution for $K$ and $M$ if and only if it has one for $cK$ and $cM$
for any $c>0$), we shall prove (i) for $\tau K$ and $\tau M$, where
\[
\tau:=\frac{2}{\left(  b-a\right)  \int_{a}^{b}K\left(  t\right)  dt}.
\]
Let $K_{a}$ and $K_{b}$ be given by (\ref{om}). Since $\delta_{\Omega}%
\leq\left(  b-a\right)  /2$ in $\Omega$ we observe that by (\ref{l11})
\[
\mathcal{S}\left(  \tau K\right)  \leq\frac{\tau}{b-a}\max\left(  \int_{a}%
^{b}K_{a}\left(  t\right)  dt,\int_{a}^{b}K_{b}\left(  t\right)  dt\right)
\delta_{\Omega}\leq\tau\delta_{\Omega}\int_{a}^{b}K\left(  t\right)  dt\leq1
\]
in $\Omega$. Let $\gamma\in\left(  0,\left(  p-1\right)  /p\right)  $, and let
us now define%
\begin{gather*}
\mathcal{M}_{p}:=\max\left(  \left\Vert M_{a}\right\Vert _{L^{p}\left(
\Omega\right)  },\left\Vert M_{b}\right\Vert _{L^{p}\left(  \Omega\right)
}\right)  \text{,\qquad}\beta:=\frac{\tau}{b-a},\\
r:=\left(  \beta\left\Vert \delta_{\Omega}^{-\gamma}\right\Vert _{L^{p^{\prime
}}\left(  \Omega\right)  }\mathcal{M}_{p}\gamma\right)  ^{1/\left(
\gamma+1\right)  }\text{,\qquad}d:=r\delta_{\Omega}\text{,}\\
\mathcal{C}:=\left\{  v\in C\left(  \overline{\Omega}\right)  :d\leq v\leq
\tau\mathcal{S}\left(  K\right)  \text{ in }\Omega\right\}  .
\end{gather*}
A simple computation shows that
\[
\left\Vert \delta_{\Omega}^{-\gamma}\right\Vert _{L^{p^{\prime}}\left(
\Omega\right)  }=2^{\gamma}\frac{\left(  b-a\right)  ^{1/p^{\prime}-\gamma}%
}{\left(  1-\gamma p^{\prime}\right)  ^{1/p^{\prime}}}%
\]
and so (\ref{m2}) says that
\begin{equation}
\frac{\left(  b-a\right)  ^{2\gamma}}{2^{\gamma}}\left\Vert \delta_{\Omega
}^{-\gamma}\right\Vert _{L^{p^{\prime}}\left(  \Omega\right)  }\mathcal{M}%
_{p}\leq\frac{\gamma^{\gamma}}{\left(  \gamma+1\right)  ^{\gamma+1}}%
\frac{\left(  \int_{a}^{b}K(t)\delta_{\Omega}(t)dt\right)  ^{1+\gamma}%
}{\left(  \int_{a}^{b}K(t)dt\right)  ^{\gamma}}.\label{qui}%
\end{equation}
Therefore, taking into account (\ref{l1}) and (\ref{qui}) we find that%
\begin{gather}
\tau\mathcal{S}\left(  K\right)  \geq\beta\left(  \int_{a}^{b}K\left(
t\right)  \delta_{\Omega}\left(  t\right)  dt\right)  \delta_{\Omega}%
\geq\nonumber\\
\beta\left(  \left\Vert \delta_{\Omega}^{-\gamma}\right\Vert _{L^{p^{\prime}%
}\left(  \Omega\right)  }\mathcal{M}_{p}\left(  \int_{a}^{b}K(t)dt\right)
^{\gamma}\frac{\left(  b-a\right)  ^{2\gamma}\left(  \gamma+1\right)
^{\gamma+1}}{2^{\gamma}\gamma^{\gamma}}\right)  ^{1/\left(  \gamma+1\right)
}\delta_{\Omega}=\nonumber\\
\left(  \beta\left\Vert \delta_{\Omega}^{-\gamma}\right\Vert _{L^{p^{\prime}%
}\left(  \Omega\right)  }\mathcal{M}_{p}\frac{\left(  \gamma+1\right)
^{\gamma+1}}{\gamma^{\gamma}}\right)  ^{1/\left(  \gamma+1\right)  }%
\delta_{\Omega}=r\frac{\gamma+1}{\gamma}\delta_{\Omega}\geq d\qquad\text{in
}\Omega\text{.}\label{beta}%
\end{gather}
In particular, $\mathcal{C}\not =\emptyset$. On the other hand, utilizing
(\ref{l11}) we also deduce that%
\begin{gather}
\tau\mathcal{S}\left(  Md^{-\gamma}\right)  \leq\beta r^{-\gamma}\max\left(
\int_{a}^{b}M_{a}\left(  t\right)  \delta_{\Omega}^{-\gamma}\left(  t\right)
dt,\int_{a}^{b}M_{b}\left(  t\right)  \delta_{\Omega}^{-\gamma}\left(
t\right)  dt\right)  \delta_{\Omega}\leq\nonumber\\
\beta r^{-\gamma}\left\Vert \delta_{\Omega}^{-\gamma}\right\Vert
_{L^{p^{\prime}}\left(  \Omega\right)  }\mathcal{M}_{p}\delta_{\Omega}%
=\frac{1}{\gamma}r\delta_{\Omega}\qquad\text{in }\Omega\text{.}%
\label{dosbetaa}%
\end{gather}
($\mathcal{S}\left(  Md^{-\gamma}\right)  $ is well defined since
$Md^{-\gamma}\in L^{s}\left(  \Omega\right)  $ for some $s>1$ because $M\in
L^{p}\left(  \Omega\right)  $ with $p\geq2$ and $\gamma p^{\prime}<1$.) For
$v\in\mathcal{C}$ we now set $u:=\tau\mathcal{S}\left(  K-Mv^{-\gamma}\right)
:=\mathcal{T}\left(  v\right)  $. Recalling (\ref{beta}), (\ref{dosbetaa}) and
that $v\geq d$ we derive that%
\[
\tau\mathcal{S}\left(  K\right)  \geq u\geq\tau\mathcal{S}\left(
K-Md^{-\gamma}\right)  \geq\left(  r\frac{\gamma+1}{\gamma}-\frac{1}{\gamma
}r\right)  \delta_{\Omega}=d\qquad\text{in }\Omega
\]
and thus $u\in\mathcal{C}$. Moreover, since $\gamma<1/p^{\prime}$ one can see
that $v\rightarrow K-Mv^{-\gamma}$ is continuous from $C\left(  \overline
{\Omega}\right)  $ into $L^{s}\left(  \Omega\right)  $ for some $s>1$, and
then it is easy to check employing the Sobolev imbedding theorems that
$\mathcal{T}:\mathcal{C}\rightarrow\mathcal{C}$ is a continuous compact
operator. Hence, Schauder's fixed point theorem gives some $v\in\mathcal{C}$
solution of
\begin{equation}
\left\{
\begin{array}
[c]{ll}%
-v^{\prime\prime}=K\left(  x\right)  -M\left(  x\right)  v^{-\gamma} &
\text{in }\Omega\\
v>0 & \text{in }\Omega\\
v=0 & \text{on }\partial\Omega.
\end{array}
\right.  \label{aux}%
\end{equation}
Furthermore, $v\in C^{1}\left(  \overline{\Omega}\right)  $, and since
$v\leq1$ (because $v\leq\mathcal{S}\left(  \tau K\right)  \leq1$) it follows
from (\ref{aux}) that $v$ is a subsolution of (\ref{prob}) for $\lambda=1$.
Therefore recalling Remark \ref{supersol} and Theorem \ref{subsup} we obtain
some $u\in C^{1}\left(  \Omega\right)  \cap C\left(  \overline{\Omega}\right)
$ solution of (\ref{prob}). Moreover, if $K\delta_{\Omega}^{-\gamma}\in
L^{r}\left(  \Omega\right)  $ for some $r>1$, since we also have that
$M\delta_{\Omega}^{-\gamma}\in L^{s}\left(  \Omega\right)  $ with $s>1$ and
$u\geq c\delta_{\Omega}$ for some $c>0$, by standard regularity arguments we
conclude that $u\in W^{2,q}\left(  \Omega\right)  \cap P^{\circ}$, $q>1$. This
ends the proof of (i).

In order to prove (ii) we proceed similarly. Without loss of generality we
assume that $\mathcal{S}\left(  K\right)  \leq1$ in $\Omega$. We now utilize
(\ref{l1}) to get that
\[
\mathcal{S}\left(  K\right)  \geq\frac{1}{b-a}\left(  \int_{a}^{b}K\left(
t\right)  \delta_{\Omega}\left(  t\right)  dt\right)  \delta_{\Omega}%
:=\frac{c_{K}}{b-a}\delta_{\Omega}\qquad\text{in }\Omega.
\]
Also, by (\ref{hip}) we may fix $\varepsilon\in\left(  0,1\right)  $ such
that
\[
\varepsilon\leq\min\left(  \frac{2}{c_{K}},\frac{2\left(  c_{K}-\mathcal{M}%
_{1}\right)  }{2c_{K}+1}\right)  ,
\]
where $\mathcal{M}_{1}$ is defined as in (i). Let $d:=\frac{\varepsilon}%
{b-a}c_{K}\delta_{\Omega}$, and note that $d\leq1$ in $\Omega$. We next choose
$\gamma_{0}>0$ such that for all $\gamma\in\left(  0,\gamma_{0}\right]  $ it
holds that $Md^{-\gamma}\in L^{s}\left(  \Omega\right)  $ with $s>1$. Making
$\gamma_{0}$ smaller if necessary, from (\ref{l11}) and Lebesgue's dominated
convergence theorem we obtain that for such $\gamma$%
\begin{gather*}
\mathcal{S}\left(  Md^{-\gamma}\right)  \leq\frac{1}{b-a}\max\left(  \int
_{a}^{b}M_{a}\left(  t\right)  d^{-\gamma}\left(  t\right)  dt,\int_{a}%
^{b}M_{b}\left(  t\right)  d^{-\gamma}\left(  t\right)  dt\right)
\delta_{\Omega}\leq\\
\frac{1}{b-a}\left(  \mathcal{M}_{1}+\frac{\varepsilon}{2}\right)
\delta_{\Omega}\leq\frac{1-\varepsilon}{b-a}c_{K}\delta_{\Omega}\qquad\text{in
}\Omega.
\end{gather*}

Define now $\mathcal{C}:=\left\{  v\in C\left(  \overline{\Omega}\right)
:d\leq v\leq\mathcal{S}\left(  K\right)  \text{ in }\Omega\right\}  $, and for
$v\in\mathcal{C}$ let $u:=\mathcal{S}\left(  K-Mv^{-\gamma}\right)  $. For
$\gamma\in\left(  0,\gamma_{0}\right]  $ we have that%
\[
\mathcal{S}\left(  K\right)  \geq u\geq\mathcal{S}\left(  K-Md^{-\gamma
}\right)  \geq\frac{1}{b-a}\left(  c_{K}-\left(  1-\varepsilon\right)
c_{K}\right)  \delta_{\Omega}\qquad\text{in }\Omega
\]
and thus $u\in\mathcal{C}$. Now the proof of (ii) can be continued exactly as
the proof of (i). $\blacksquare$

\begin{remark}
\label{compar}Let us mention that one can verify that the inequality
(\ref{hip}) is better than (\ref{m2}), but on the other hand in (ii) there is
no lower estimate for $\gamma_{0}$. Let us also note that if for instance
$\gamma<1/2$ then $K\delta_{\Omega}^{-\gamma}\in L^{r}\left(  \Omega\right)  $
for some $r>1$ (in fact, if $K\in L^{p}\left(  \Omega\right)  $ with $p\geq2$,
then $K\delta_{\Omega}^{-\gamma}\in L^{r}\left(  \Omega\right)  $ with $r>1$
when $p\left(  1-\gamma\right)  >1$). $\blacksquare$
\end{remark}

\begin{corollary}
\label{coro}Let $K$, $M$ and $\gamma$ be as in Theorem \ref{nuevo}, and let
$\alpha>\gamma$. Then there exists $\lambda_{0}>0$ such that the problem
(\ref{prob}) has a solution $u\in C^{1}\left(  \Omega\right)  \cap C\left(
\overline{\Omega}\right)  $ for all $\lambda\leq\lambda_{0}$; and $u\in
W^{2,q}\left(  \Omega\right)  \cap P^{\circ}$, $q>1$, whenever $K\delta
_{\Omega}^{-\alpha}\in L^{r}\left(  \Omega\right)  $ with $r>1$.
\end{corollary}

\textit{Proof}. As in the above theorem it suffices to construct a subsolution
for (\ref{prob}). Let $u$ be the solution of (\ref{prob}) for $\lambda=1$ and
$\alpha=\gamma$ provided by either Theorem \ref{nuevo} (i) or (ii). We choose
$0<\varepsilon\leq\min\left(  1,1/\left\Vert u\right\Vert _{\infty}\right)  $.
Now, for every $\alpha>\gamma$ we get that
\begin{gather*}
-\left(  \varepsilon u\right)  ^{\prime\prime}=\left(  \varepsilon^{1+\gamma
}K\left(  x\right)  -\varepsilon^{1+\gamma}M\left(  x\right)  \right)  \left(
\varepsilon u\right)  ^{-\gamma}\leq\\
K\left(  x\right)  \left(  \varepsilon u\right)  ^{-\alpha}-\varepsilon
^{1+\gamma}M\left(  x\right)  \left(  \varepsilon u\right)  ^{-\gamma}%
\qquad\text{in }\Omega
\end{gather*}
and hence $\varepsilon u$ is a subsolution of (\ref{prob}) for $\lambda
_{0}:=\varepsilon^{1+\gamma}$. $\blacksquare$

\qquad

In the next theorem we complement the results contained in Theorem \ref{nuevo}
(ii) and Corollary \ref{coro}, without imposing any relation between $\alpha$
and $\gamma$. We set
\[
\mathcal{M}:=\left\{  x\in\Omega:M\left(  x\right)  >0\right\}  .
\]

\begin{theorem}
\label{algo}Let $\alpha,\gamma>0$, let $0\leq K,M\in L^{2}\left(
\Omega\right)  $ with $M\in C\left(  \Omega\right)  $ and let $K_{a}$, $K_{b}%
$, $M_{a}$, $M_{b}$ be given by (\ref{om}). Suppose that $\overline
{\mathcal{M}}\subset\Omega$ and that (\ref{hip}) holds. Then there exists
\begin{equation}
\lambda_{0}\geq\left(  \frac{dist\left(  \mathcal{M},\partial\Omega\right)
}{b-a}\right)  ^{\gamma}\frac{\left(  \int_{a}^{b}K\delta_{\Omega}-\max\left(
\int_{a}^{b}M_{a},\int_{a}^{b}M_{b}\right)  \right)  ^{\gamma}}{\left(
\frac{1}{2}\left(  \max\left(  \int_{a}^{b}K_{a},\int_{a}^{b}K_{b}\right)
-\int_{a}^{b}M\delta_{\Omega}\right)  \right)  ^{\alpha\frac{1+\gamma
}{1+\alpha}}} \label{lambda}%
\end{equation}
such that for all $\lambda\leq\lambda_{0}$ the problem (\ref{prob}) has a
solution $u\in C^{1}\left(  \Omega\right)  \cap C\left(  \overline{\Omega
}\right)  $; and $u\in W^{2,q}\left(  \Omega\right)  \cap P^{\circ}$, $q>1$,
whenever $K\delta_{\Omega}^{-\alpha},M\delta_{\Omega}^{-\gamma}\in
L^{r}\left(  \Omega\right)  $ with $r>1$.
\end{theorem}

\textit{Proof}. Let $\phi:=\mathcal{S}\left(  K-M\right)  $. Applying Lemma
\ref{tomas} to $\mathcal{S}\left(  K\right)  $ and $\mathcal{S}\left(
M\right)  $ we see that in $\Omega$%
\begin{gather*}
\phi\geq\frac{1}{b-a}\left(  \int_{a}^{b}K\left(  t\right)  \delta_{\Omega
}\left(  t\right)  dt-\max\left(  \int_{a}^{b}M_{a}\left(  t\right)
dt,\int_{a}^{b}M_{b}\left(  t\right)  dt\right)  \right)  \delta_{\Omega},\\
\phi\leq\frac{1}{b-a}\left(  \max\left(  \int_{a}^{b}K_{a}\left(  t\right)
dt,\int_{a}^{b}K_{b}\left(  t\right)  dt\right)  -\int_{a}^{b}M\left(
t\right)  \delta_{\Omega}\left(  t\right)  dt\right)  \delta_{\Omega}.
\end{gather*}
We note that in particular it follows from (\ref{hip}) that $\phi>0$ in
$\Omega$. Let us now define
\begin{gather*}
\overline{\mu}:=\left(  \frac{1}{2}\left(  \max\left(  \int_{a}^{b}K_{a}%
,\int_{a}^{b}K_{b}\right)  -\int_{a}^{b}M\delta_{\Omega}\right)  \right)
^{\alpha},\\
\underline{\beta}:=\left(  \frac{dist\left(  \mathcal{M},\partial
\Omega\right)  }{b-a}\left(  \int_{a}^{b}K\delta_{\Omega}-\max\left(  \int
_{a}^{b}M_{a},\int_{a}^{b}M_{b}\right)  \right)  \right)  ^{\gamma}.
\end{gather*}
Let $\mu\geq\overline{\mu}$ and $\beta\leq\underline{\beta}$ and set
$\varepsilon:=\mu^{-1/\left(  1+\alpha\right)  }$. Since $\phi\leq
\overline{\mu}^{1/\alpha}$ in $\Omega$ and $\phi\geq\underline{\beta
}^{1/\gamma}$ in $\mathcal{M}$ we have%
\begin{gather*}
-\left(  \varepsilon\phi\right)  ^{\prime\prime}=\mu^{-1}\varepsilon^{-\alpha
}K\left(  x\right)  -\varepsilon^{1+\gamma}\varepsilon^{-\gamma}M\left(
x\right)  \leq\\
K\left(  x\right)  \left(  \varepsilon\phi\right)  ^{-\alpha}-\left(
\frac{\beta}{\mu^{\left(  1+\gamma\right)  /\left(  1+\alpha\right)  }%
}\right)  M\left(  x\right)  \left(  \varepsilon\phi\right)  ^{-\gamma}%
\qquad\text{in }\Omega
\end{gather*}
and so $\varepsilon\phi$ is a subsolution of (\ref{prob}) for $\lambda
=\beta\mu^{-\left(  1+\gamma\right)  /\left(  1+\alpha\right)  }$. Since
(\ref{lambda}) is clearly true the theorem follows. $\blacksquare$

\qquad

The next theorem states necessary conditions on $\lambda$, $K$ and $M$ in
order for (\ref{prob}) to admit solutions in the case $\alpha\leq\gamma$.

\begin{theorem}
\label{necesaria}Let $0<\alpha\leq\gamma$ and let $K_{a}$ and $K_{b}$ be given
by (\ref{om}). Suppose for some $\lambda>0$ (\ref{prob}) has a solution $u\in
W^{2,q}\left(  \Omega\right)  $, $q>1$. Then
\[
\lambda<\left(  \frac{\left(  \alpha+1\right)  }{2}\max\left(  \int_{a}%
^{b}K_{a}\left(  t\right)  dt,\int_{a}^{b}K_{b}\left(  t\right)  dt\right)
\right)  ^{\frac{\gamma-\alpha}{\alpha+1}}\frac{\int_{a}^{b}K\left(  t\right)
\delta_{\Omega}\left(  t\right)  dt}{\int_{a}^{b}M\left(  t\right)
\delta_{\Omega}\left(  t\right)  dt}.
\]

\end{theorem}

\textit{Proof}. Let $u>0$ be a solution of (\ref{prob}) for some $\lambda>0$
and pick $\sigma:=\alpha+1$. We get that
\begin{gather*}
-\left(  u^{\sigma}\right)  ^{\prime\prime}=-\sigma u^{\sigma-1}%
u^{\prime\prime}-\sigma\left(  \sigma-1\right)  u^{\sigma-2}\left(  u^{\prime
}\right)  ^{2}\leq\\
-\sigma u^{\sigma-1}u^{\prime\prime}\leq\sigma u^{\sigma-1}K\left(  x\right)
u^{-\alpha}=\sigma K\left(  x\right)  \qquad\text{in }\Omega
\end{gather*}
and so, since $\mathcal{S}$ is a positive operator, recalling (\ref{l11}) we
find that
\begin{gather*}
0<u^{\sigma}\leq\sigma\mathcal{S}\left(  K\right)  \leq\frac{\sigma}{b-a}%
\max\left(  \int_{a}^{b}K_{a}\left(  t\right)  dt,\int_{a}^{b}K_{b}\left(
t\right)  dt\right)  \delta_{\Omega}\leq\\
\frac{\sigma}{2}\max\left(  \int_{a}^{b}K_{a}\left(  t\right)  dt,\int_{a}%
^{b}K_{b}\left(  t\right)  dt\right)  \qquad\text{in }\Omega.
\end{gather*}
Therefore, it follows that
\[
\left\Vert u\right\Vert _{\infty}\leq\left(  \frac{\left(  \alpha+1\right)
}{2}\max\left(  \int_{a}^{b}K_{a}\left(  t\right)  dt,\int_{a}^{b}K_{b}\left(
t\right)  dt\right)  \right)  ^{1/\left(  \alpha+1\right)  }.
\]

Let $\varepsilon:=1/\left\Vert u\right\Vert _{\infty}$. On the other side,
since $\alpha\leq\gamma$ and $\sigma-1=\alpha$ we also have that
\begin{gather*}
-\left(  \left(  \varepsilon u\right)  ^{\sigma}\right)  ^{\prime\prime}%
\leq-\sigma\left(  \varepsilon u\right)  ^{\sigma-1}\left(  \varepsilon
u\right)  ^{\prime\prime}\leq\\
\sigma\left(  \varepsilon u\right)  ^{\sigma-1}\left(  \varepsilon^{1+\alpha
}K\left(  x\right)  -\lambda\varepsilon^{1+\gamma}M\left(  x\right)  \right)
\left(  \varepsilon u\right)  ^{-\alpha}=\\
\sigma\left(  \varepsilon^{1+\alpha}K\left(  x\right)  -\lambda\varepsilon
^{1+\gamma}M\left(  x\right)  \right)  \qquad\text{in }\Omega.
\end{gather*}
Hence the positivity of $\mathcal{S}$ now tells us that $\mathcal{S}\left(
\varepsilon^{1+\alpha}K-\lambda\varepsilon^{1+\gamma}M\right)  >0$ and thus
$\lambda<\varepsilon^{-\left(  \gamma-\alpha\right)  }\mathcal{S}\left(
K\right)  /\mathcal{S}\left(  M\right)  $. Furthermore, since
\[
\mathcal{S}\left(  K\right)  \left(  \frac{a+b}{2}\right)  =\frac{1}{2}%
\int_{a}^{b}K\left(  t\right)  \delta_{\Omega}\left(  t\right)  dt
\]
(see the proof of Lemma \ref{tomas}) and an analogous statement is valid for
$\mathcal{S}\left(  M\right)  $, the theorem follows employing the upper bound
for $\left\Vert u\right\Vert _{\infty}$ derived in the first part of the
proof. $\blacksquare$

\begin{remark}
\label{rem}Given $M$, $K$, $\alpha$ and $\gamma$, let $\Lambda:=\left\{
\lambda>0:\left(  \ref{prob}\right)  \text{ has a solution}\right\}  $. If
$\Lambda\not =\emptyset$ (for instance if (\ref{hip}) holds) then either
$\Lambda=\left(  0,\lambda_{0}\right)  $ or $\left(  0,\lambda_{0}\right]  $
for some $0<\lambda_{0}\leq\infty$. Indeed, define $\lambda_{0}:=\sup
_{\lambda\in\Lambda}$, and now this follows from Theorem \ref{subsup}, Remark
\ref{supersol} and the fact any solution of (\ref{prob}) is a subsolution of
(\ref{prob}) with $\underline{\lambda}$ in place of $\lambda\,$whenever
$\underline{\lambda}\leq\lambda$. Let us note that if $\alpha\leq\gamma$ the
above theorem says that $\lambda_{0}<\infty$. $\blacksquare$
\end{remark}

Let now $m:\Omega\rightarrow\mathbb{R}$ be a function that may change sign in
$\Omega$. We write as usual $m=m^{+}-m^{-}$ with $m^{+}:=\max\left(
m,0\right)  $ and $m^{-}:=\max\left(  -m,0\right)  $. As a direct consequence
of Theorems \ref{nuevo} and \ref{necesaria} we obtain the

\begin{corollary}
\label{golazo}Let $0\not \equiv m\in L^{2}\left(  \Omega\right)  $ and let
$m_{a}^{-}$ and $m_{b}^{-}$ be given by (\ref{om}).\newline(i) If $m^{-}\in
L^{p}\left(  \Omega\right)  $ with $p\geq2$, $\gamma\in\left(  0,\left(
p-1\right)  /p\right)  $ and%
\begin{gather*}
\max\left(  \left\Vert m_{a}^{-}\right\Vert _{L^{p}\left(  \Omega\right)
},\left\Vert m_{b}^{-}\right\Vert _{L^{p}\left(  \Omega\right)  }\right)  \leq
c_{\gamma,p,a,b}\frac{\left(  \int_{a}^{b}m^{+}(t)\delta_{\Omega}(t)dt\right)
^{1+\gamma}}{\left(  \int_{a}^{b}m^{+}(t)dt\right)  ^{\gamma}},\text{\qquad
where}\\
c_{\gamma,p,a,b}:=\frac{\gamma^{\gamma}}{\left(  \gamma+1\right)  ^{\gamma+1}%
}\frac{\left(  1-\gamma p^{\prime}\right)  ^{1/p^{\prime}}}{\left(
b-a\right)  ^{\gamma+1/p^{\prime}}},
\end{gather*}
then (\ref{sipi}) has a solution $u\in C^{1}\left(  \Omega\right)  \cap
C\left(  \overline{\Omega}\right)  $; and $u\in W^{2,q}\left(  \Omega\right)
\cap P^{\circ}$, $q>1$, whenever $K\delta_{\Omega}^{-\gamma}\in L^{r}\left(
\Omega\right)  $ with $r>1$.\newline(ii) If
\[
\max\left(  \int_{a}^{b}m_{a}^{-}(t)dt,\int_{a}^{b}m_{b}^{-}(t)dt\right)
<\int_{a}^{b}m^{+}\left(  t\right)  \delta_{\Omega}\left(  t\right)  dt,
\]
then (\ref{sipi}) has a solution $u\in W^{2,q}\left(  \Omega\right)  \cap
P^{\circ}$, $q>1$, for all $\gamma\in\left(  0,\gamma_{0}\right]  $ and some
$\gamma_{0}>0$. \newline(iii) If (\ref{sipi}) has a solution $u\in
W^{2,q}\left(  \Omega\right)  $, $q>1$, then
\[
\int_{a}^{b}m^{-}\left(  t\right)  \delta_{\Omega}\left(  t\right)
dt<\int_{a}^{b}m^{+}\left(  t\right)  \delta_{\Omega}\left(  t\right)  dt.
\]

\end{corollary}

\section{The $N$-dimensional problem}

We consider now the case of a smooth bounded domain $\mathbb{R}^{N}$, $N\geq
2$. We shall denote
\begin{gather*}
B_{R}:=\left\{  x\in\mathbb{R}^{N}:\left\vert x\right\vert <R\right\}  ,\\
\omega_{N-1}:=\text{surface area of the unit sphere }\partial B_{1}\text{ in
}\mathbb{R}^{N},\\
diam\left(  \Omega\right)  :=\sup_{x,y\in\Omega}\left\vert x-y\right\vert ,\\
P^{\circ}:=\text{interior of the positive cone of }C^{1}\left(  \overline
{\Omega}\right)  .
\end{gather*}
Since all the proofs in this section are natural adaptations of the proofs in
the one-dimensional case, we shall only give the minor changes that are
needed. For simplicity in the first theorem we shall assume that $M\in
L^{\infty}(\Omega)$.

\begin{theorem}
\label{bola}Let $0\leq K\in L^{p}\left(  \Omega\right)  $, $p>N$, and let
$0\leq M\in L^{\infty}\left(  \Omega\right)  $. Assume $\alpha=\gamma$ and let
$c_{\Omega}$ and $C_{\Omega}$ be given by Lemma \ref{morel}. If $\gamma
\in\left(  0,1/N\right)  $ and%
\begin{gather}
\left\Vert M\right\Vert _{L^{\infty}\left(  \Omega\right)  }<c_{\Omega
,\gamma,N}\frac{\left(  \int_{\Omega}K\left(  x\right)  \delta_{\Omega}\left(
x\right)  dx\right)  ^{1+\gamma}}{\left\Vert K\right\Vert _{L^{p}\left(
\Omega\right)  }^{\gamma}},\text{\qquad where}\label{hipp}\\
c_{\Omega,\gamma,N}:=\left(  \frac{c_{\Omega}}{C_{\Omega}}\right)  ^{1+\gamma
}\frac{\gamma^{\gamma}}{\left(  \gamma+1\right)  ^{\gamma+1}}\left(  \frac
{2}{diam\left(  \Omega\right)  }\right)  ^{\gamma}\frac{1}{\left\Vert
\delta_{\Omega}^{-\gamma}\right\Vert _{L^{N}\left(  \Omega\right)  }%
},\nonumber
\end{gather}
then (\ref{prob}) has a solution $u\in C^{1}\left(  \Omega\right)  \cap
C\left(  \overline{\Omega}\right)  $ for all $\lambda\leq1$; and $u\in
W^{2,q}\left(  \Omega\right)  \cap P^{\circ}$, $q>N$, whenever $K\delta
_{\Omega}^{-\gamma}\in L^{r}\left(  \Omega\right)  $ with $r>N$. In
particular, if $\Omega:=B_{R}$ then%
\begin{equation}
c_{\Omega,\gamma,N}\geq\left(  \frac{c_{\Omega}}{C_{\Omega}}\right)
^{1+\gamma}\frac{\gamma^{\gamma}}{\left(  \gamma+1\right)  ^{\gamma+1}}%
\frac{1}{R}\left(  \frac{1-\gamma N}{\omega_{N-1}}\right)  ^{1/N}.
\label{hippp}%
\end{equation}

\end{theorem}

\textit{Proof}. For $\gamma\in\left(  0,1/N\right)  $ we pick $q\in\left(
N,1/\gamma\right)  $ and we set%
\begin{gather*}
\tau:=\frac{2}{C_{\Omega}\left\Vert K\right\Vert _{L^{p}\left(  \Omega\right)
}diam\left(  \Omega\right)  }\text{,\qquad}\mu:=\tau C_{\Omega}\left\Vert
M\right\Vert _{L^{\infty}\left(  \Omega\right)  }\left\Vert \delta_{\Omega
}^{-\gamma}\right\Vert _{L^{q}\left(  \Omega\right)  },\\
\beta:=\tau c_{\Omega}\int_{\Omega}K\left(  x\right)  \delta_{\Omega}\left(
x\right)  dx\text{,\qquad}r:=\left(  \mu\gamma\right)  ^{1/\left(
\gamma+1\right)  }\text{,\qquad}d:=r\delta_{\Omega}\text{,}\\
\mathcal{C}:=\left\{  v\in C\left(  \overline{\Omega}\right)  :d\leq v\leq
\tau\mathcal{S}\left(  K\right)  \text{ in }\Omega\right\}  .
\end{gather*}
Taking into account Lemma \ref{morel} we find that $\mathcal{S}\left(  \tau
K\right)  \leq1$ in $\Omega$ and also%
\begin{gather*}
\mathcal{S}\left(  \tau K\right)  \geq\beta\delta_{\Omega}\qquad\text{and}\\
\mathcal{S}\left(  \tau Md^{-\gamma}\right)  \leq\tau\left\Vert M\right\Vert
_{L^{\infty}\left(  \Omega\right)  }\mathcal{S}\left(  d^{-\gamma}\right)
\leq\mu r^{-\gamma}\delta_{\Omega}\qquad\text{in }\Omega.
\end{gather*}
On the other hand, we may fix $q$ close enough to $N$ so that (\ref{hipp})
implies that
\[
\mu\leq\beta^{1+\gamma}\frac{\gamma^{\gamma}}{\left(  \gamma+1\right)
^{\gamma+1}}.
\]
For $v\in\mathcal{C}$ define now $u:=\tau\mathcal{S}\left(  K-Mv^{-\gamma
}\right)  $. Taking into account the aforementioned facts we derive that%
\begin{gather*}
\tau\mathcal{S}\left(  K\right)  \geq u\geq\tau\mathcal{S}\left(
K-Md^{-\gamma}\right)  \geq\left(  \beta-\mu r^{-\gamma}\right)
\delta_{\Omega}=\\
\left(  \mu\gamma\right)  ^{1/\left(  \gamma+1\right)  }\left(  \frac{\beta
}{\left(  \mu\gamma\right)  ^{1/\left(  \gamma+1\right)  }}-\frac{1}{\gamma
}\right)  \delta_{\Omega}\geq\left(  \mu\gamma\right)  ^{1/\left(
\gamma+1\right)  }\left(  \frac{\gamma+1}{\gamma}-\frac{1}{\gamma}\right)
\delta_{\Omega}=d
\end{gather*}
in $\Omega$ and hence $u\in\mathcal{C}$. Furthermore, since $\gamma<1/N$ and
$v\geq d$ it holds that $v\rightarrow K-Mv^{-\gamma}$ is continuous from
$C\left(  \overline{\Omega}\right)  $ into $L^{s}\left(  \Omega\right)  $ for
some $s>N$, and then the proof can be continued as the proof of Theorem
\ref{nuevo} (i).

Finally, if $\Omega:=B_{R}$, using polar coordinates one can see that%
\begin{gather*}
\left\Vert \delta_{\Omega}^{-\gamma}\right\Vert _{L^{N}\left(  \Omega\right)
}^{N}=\omega_{N-1}\int_{0}^{R}r^{N-1}\left(  R-r\right)  ^{-\gamma N}dr\leq\\
\omega_{N-1}R^{N-1}\int_{0}^{R}r^{-\gamma N}dr=\frac{\omega_{N-1}}{1-\gamma
N}R^{N\left(  1-\gamma\right)  }%
\end{gather*}
which in turn yields (\ref{hippp}). $\blacksquare$

\qquad

With the same proof as in the previous section we deduce the

\begin{corollary}
\label{coro n}Let $K$, $M$ and $\gamma$ be as in Theorem \ref{bola}, and let
$\alpha>\gamma$. Then there exists $\lambda_{0}>0$ such that the problem
(\ref{prob}) has a solution $u\in C^{1}\left(  \Omega\right)  \cap C\left(
\overline{\Omega}\right)  $ for all $\lambda\leq\lambda_{0}$; and $u\in
W^{2,q}\left(  \Omega\right)  \cap P^{\circ}$, $q>N$, whenever $K\delta
_{\Omega}^{-\alpha}\in L^{r}\left(  \Omega\right)  $ with $r>N$.
\end{corollary}

Let $\mathcal{M}:=\left\{  x\in\Omega:M\left(  x\right)  >0\right\}  $. We
also have

\begin{theorem}
\label{algo n}Let $\alpha,\gamma>0$ and let $0\leq K,M\in L^{p}\left(
\Omega\right)  $, $p>N$, and $M\in C\left(  \Omega\right)  $. Let $c_{\Omega}$
and $C_{\Omega}$ be given by Lemma \ref{morel}. Suppose that $\overline
{\mathcal{M}}\subset\Omega$ and that
\[
\left\Vert M\right\Vert _{L^{p}\left(  \Omega\right)  }<\frac{c_{\Omega}%
}{C_{\Omega}}\int_{\Omega}K\left(  x\right)  \delta_{\Omega}\left(  x\right)
dx.
\]
Then there exists
\[
\lambda_{0}\geq\left(  dist\left(  \mathcal{M},\partial\Omega\right)  \right)
^{\gamma}\frac{\left(  c_{\Omega}\int_{\Omega}K\delta_{\Omega}-C_{\Omega
}\left\Vert M\right\Vert _{L^{p}\left(  \Omega\right)  }^{\gamma}\right)
^{\gamma}}{\left(  \frac{diam\left(  \Omega\right)  }{2}\left(  C_{\Omega
}\left\Vert K\right\Vert _{L^{p}\left(  \Omega\right)  }-c_{\Omega}%
\int_{\Omega}M\delta_{\Omega}\right)  \right)  ^{\alpha\frac{1+\gamma
}{1+\alpha}}}%
\]
such that for all $\lambda\leq\lambda_{0}$ the problem (\ref{prob}) has a
solution $u\in C^{1}\left(  \Omega\right)  \cap C\left(  \overline{\Omega
}\right)  $; and $u\in W^{2,q}\left(  \Omega\right)  \cap P^{\circ}$, $q>N$,
whenever $K\delta_{\Omega}^{-\alpha},M\delta_{\Omega}^{-\gamma}\in
L^{r}\left(  \Omega\right)  $ with $r>N$.
\end{theorem}

\textit{Proof}. The proof follows as in Theorem \ref{algo} defining now
\begin{gather*}
\overline{\mu}:=\left(  \frac{diam\left(  \Omega\right)  }{2}\left(
C_{\Omega}\left\Vert K\right\Vert _{L^{p}\left(  \Omega\right)  }-c_{\Omega
}\int_{\Omega}M\delta_{\Omega}\right)  \right)  ^{\alpha}\\
\underline{\beta}:=\left(  dist\left(  \mathcal{M},\partial\Omega\right)
\left(  c_{\Omega}\int_{\Omega}K\delta_{\Omega}-C_{\Omega}\left\Vert
M\right\Vert _{L^{p}\left(  \Omega\right)  }\right)  \right)  ^{\gamma}%
.\qquad\blacksquare
\end{gather*}

\begin{theorem}
\label{necesaria n}Let $0<\alpha\leq\gamma$ and let $C_{\Omega}$ be given by
Lemma \ref{morel}. Suppose for some $\lambda>0$ (\ref{prob}) has a solution
$u\in W^{2,q}\left(  \Omega\right)  $, $q>N$. Then
\begin{equation}
\lambda<\left(  \frac{diam\left(  \Omega\right)  }{2}C_{\Omega}\left(
\alpha+1\right)  \left\Vert K\right\Vert _{L^{p}\left(  \Omega\right)
}\right)  ^{\frac{\gamma-\alpha}{\alpha+1}}\inf_{\Omega}\left(  \frac
{\mathcal{S}\left(  K\right)  }{\mathcal{S}\left(  M\right)  }\right)  .
\label{n}%
\end{equation}

\end{theorem}

\textit{Proof}. Suppose $u$ is a solution of (\ref{prob}) for some $\lambda
>0$. Employing Lemma \ref{morel} and the positivity of $\mathcal{S}$ and
arguing as in Theorem \ref{necesaria} we can prove that
\begin{gather*}
\left\Vert u\right\Vert _{\infty}\leq\left(  \frac{diam\left(  \Omega\right)
}{2}C_{\Omega}\left(  \alpha+1\right)  \left\Vert K\right\Vert _{L^{p}\left(
\Omega\right)  }\right)  ^{1/\left(  \alpha+1\right)  }\qquad\text{and}\\
\lambda<\left\Vert u\right\Vert _{\infty}^{\gamma-\alpha}\inf_{\Omega}\left(
\frac{\mathcal{S}\left(  K\right)  }{\mathcal{S}\left(  M\right)  }\right)
\end{gather*}
and this gives (\ref{n}). $\blacksquare$

\begin{remark}
Let us note that the statement in Remark \ref{rem} is clearly also valid for
the $N$-dimensional problem. $\blacksquare$
\end{remark}

\begin{corollary}
\label{golazo n}Let $m\in L^{p}\left(  \Omega\right)  $ with $p>N$ and
$m^{-}\in L^{\infty}\left(  \Omega\right)  $. Let $c_{\Omega}$ and $C_{\Omega
}$ be given by Lemma \ref{morel}.\newline(i) If $\gamma\in\left(
0,1/N\right)  $ and%
\begin{gather*}
\left\Vert m^{-}\right\Vert _{L^{\infty}\left(  \Omega\right)  }%
<c_{\Omega,\gamma,N}\frac{\left(  \int_{\Omega}m^{+}\left(  x\right)
\delta_{\Omega}\left(  x\right)  dx\right)  ^{1+\gamma}}{\left\Vert
m^{+}\right\Vert _{L^{p}\left(  \Omega\right)  }^{\gamma}},\text{\qquad
where}\\
c_{\Omega,\gamma,N}:=\left(  \frac{c_{\Omega}}{C_{\Omega}}\right)  ^{1+\gamma
}\frac{\gamma^{\gamma}}{\left(  \gamma+1\right)  ^{\gamma+1}}\left(  \frac
{2}{diam\left(  \Omega\right)  }\right)  ^{\gamma}\frac{1}{\left\Vert
\delta_{\Omega}^{-\gamma}\right\Vert _{L^{N}\left(  \Omega\right)  }},
\end{gather*}
then (\ref{sipi}) has a solution $u\in C^{1}\left(  \Omega\right)  \cap
C\left(  \overline{\Omega}\right)  $; and $u\in W^{2,q}\left(  \Omega\right)
\cap P^{\circ}$, $q>N$, whenever $m^{+}\delta_{\Omega}^{-\gamma}\in
L^{r}\left(  \Omega\right)  $ with $r>N$. In particular, if $\Omega:=B_{R}$
then%
\[
c_{\Omega,\gamma,N}\geq\left(  \frac{c_{\Omega}}{C_{\Omega}}\right)
^{1+\gamma}\frac{\gamma^{\gamma}}{\left(  \gamma+1\right)  ^{\gamma+1}}%
\frac{1}{R}\left(  \frac{1-\gamma N}{\omega_{N-1}}\right)  ^{1/N}.
\]
\newline(ii) If (\ref{sipi}) has a solution $u\in W^{2,q}\left(
\Omega\right)  $ with $q>N$, then $\mathcal{S}\left(  m\right)  >0$ in
$\Omega$.
\end{corollary}

\end{document}